\documentclass[11pt]{amsart}
\usepackage{color}
\usepackage{graphicx}
\usepackage{hyperref}
\usepackage{amssymb}
\usepackage{amsmath}
\usepackage{amsthm}
\usepackage{mathtools}
\usepackage{fullpage}
\allowdisplaybreaks[4]

\theoremstyle{plain}
\newtheorem{theorem}{Theorem}
\newtheorem{lemma}[theorem]{Lemma}

\theoremstyle{definition}
\newtheorem*{definition*}{Definition}

\newtheorem*{notation*}{Notation}
\theoremstyle{remark}
\newtheorem*{remark*}{Remark}

\makeatletter
\@namedef{subjclassname@2010}{%
  \textup{2010} Mathematics Subject Classification}
  \makeatother

\newcommand{\R}{\mathbb{R}} 
 \newcommand{\N}{\mathbb{N}}

\begin{document}

\title{A Bound for Stieltjes constants}

\author{S. Pauli}
\author{F. Saidak}

\address{Department of Mathematics and Statistics,
   UNC Greensboro, Greensboro, NC 27402, USA }.
\email{s\_pauli@uncg.edu, f\_saidak@uncg.edu}


\begin{abstract}
The main goal of this note is to improve the best known bounds for the Stieltjes constants,
using the method of steepest descent that was applied in 2011 by Coffey and Knessl in order to approximate these constants. 
\end{abstract}

\maketitle

\section{Introduction}

Let $\zeta(s)$ denote the Riemann zeta function (for $s \in \mathbb{C}$ defined by Riemann \cite{r-0} in 1859). 
The function $\zeta(s)$ has an Euler product (Euler \cite{e-1737} of 1737) and also satisfies a 
functional equation (Euler \cite{e-1749} of 1749). In this paper we consider the related Hurwitz zeta function 
$\zeta(s, a)$ (see Hurwitz \cite{h-1882} of 1882), which for $0 < a \leq 1$ has the Laurent series expansion:
$$ \zeta(s,a) := \sum_{n=0}^{\infty} \frac{1}{(n+a)^s} = \frac{1}{s-1} + \sum_{n=0}^{\infty} \frac{(-1)^n \gamma_n(a)}{n!} (s-1)^n, $$
where, for non-negative integers $n$, the coefficients $\gamma_n(a)$ are known as the Stieltjes constants (\cite{s1}), 
which were generalized to {\em fractional} values $\alpha \in \mathbb{R}^+$ by Kreminski \cite{kreminski03} in 2003. 
These constants have several interesting and unexpected applications in the zeta function theory, as was shown recently in 
\cite{ck}, \cite{c1}, \cite{fps}, and \cite{fpsone}. Moreover, the classical Euler-Maclaurin Summation can be used to prove 
(see \cite{fpsbound}) that, if we set \(C_\alpha(a) = \gamma_\alpha(a) -  \frac{\log^\alpha(a)}{a}\) and let 
$f_{\alpha}(x) = \frac{\log^{\alpha}(x+a)}{x+a}$, then we have:
\begin{align*}
C_\alpha(a) \; =  & \; \sum\limits_{r=1}^m \frac{\log^{\alpha}(r+a)}{r+a} \;  -  \; 
\frac{\log^{\alpha+1}(m+a)}{\alpha+1} \; - \; \frac{\log^{\alpha}(m+a)}{2(m+a)}   \\
& - \; \sum\limits_{j=1}^{\lfloor v/2\rfloor} \frac{B_{2j}}{(2j)!}f_\alpha^{(2j-1)}(m)
\; + \; (-1)^{v-1}\int\limits_m^\infty P_v(x)f_\alpha^{(v)}(x) \; dx,
\end{align*}
where the $B_j$ denote the Bernoulli numbers (introduced by Bernoulli in \cite{b-1713} of 1713), and \(P_v\) is the \(v\)-th periodic Bernoulli function
(see \cite{k-1935}). This expression has many useful applications; in our recent work, we have used it to find zero-free regions for the fractional 
derivatives of the Riemann zeta function. There one of the key estimates (Lemma 4.1, \cite{fpsone}) was the bound, for $0 < \alpha \leq 1$, 
\[
\left| \int_1^{\infty}   P_3(x) f_{\alpha}''' (x) \; dx \right| <  0.013.
\]
Here, with different goals in mind, we will consider another special case of the above Euler-Maclaurin summation formula.  We set \(m=1\) and \(v=2\) and analyze the expression
\begin{align}\label{eq:Calpha}
C_\alpha(a) \; = \; & \; \frac{\log^{\alpha}(1+a)}{1+a}  \; -
 \; \frac{\log^{\alpha+1}(1+a)}{\alpha+1} \; - \; \frac{\log^{\alpha}(1+a)}{2(1+a)} \; - \; \frac{1}{12}f_\alpha'(1) \; - \; \int\limits_1^\infty P_2(x)f_\alpha''(x) \; dx.
\end{align}
Now, bounding the generalized fractional Stieltjes constants $\gamma_\alpha(a)$ (or the functions $C_\alpha(a)$), 
means finding  (this time for $1 \leq \alpha \in \mathbb{R}$) effective bounds for:
\begin{equation}\label{eq:P2}
 \left| \int_1^{\infty}   P_2(x) f_{\alpha}'' (x) \; dx \right|.
\end{equation}
Since the Bernoulli periodic function $P_2(x)$ involved in the integrals (\ref{eq:Calpha}) and (\ref{eq:P2}) has a simple Fourier series expansion, 
established by Hurwitz in 1890 (see \cite{leh:1988}) we have 
\[P_2(x) = 
\frac{-n!}{(2\pi i)^n}\sum_{k = - \infty \atop{k \neq 0}}^{\infty}\frac{e^{2\pi ikx}}{k^2} \; \Bigg\rvert_{n=2}
=
\frac{2}{(2\pi)^2}\sum_{k=1}^\infty\frac{e^{2\pi ikx}+e^{-2\pi ikx}}{k^2}
=
\frac{4}{(2\pi)^2}\sum_{k=1}^\infty\frac{\Re(e^{2\pi ikx})}{k^2}.
\]
Because this series is absolutely convergent and \(f''_\alpha(x)\) is bounded (see (\ref{eq ddf})), if we set  
\(S_k :=\int_1^\infty  e^{2\pi ikx} f_{\alpha}''(x) dx\) and $S^* := \sup_{k \in \mathbb{N}} |S_k|$, then we obtain
\begin{equation}
\label{eq:S}
\left|\int_1^\infty P_2(x) f_{\alpha}''(x) \; dx\right|  
= \left|\frac{1}{\pi^2}\sum_{k=1}^\infty\frac{1}{k^2} S_k \right|
\le \frac{1}{\pi^2}\sum_{k=1}^\infty\frac{1}{k^2} \left| S_k \right| \leq \frac{1}{6} S^*. 
\end{equation}

We obtain effective bounds on $|S_k|$ and $S^*$ by choosing a suitable integration path as done in the method of steepest decent 
for approximating integrals, see \cite[Section II 4]{wong} for example.
This path originates at a point \(b\) on the real axis and then joins a level curve with constant imaginary part
that crosses a saddle point and on the right half plane has the asymptote \(\frac{\pi}{2} i\).

We find the location of this important saddle point using the Lambert $W$ function (Lambert \cite{lam-1758} of 1758), which is 
the solution of the  special case $x = W(x) e^{W(x)}$ of the so-called Lambert transcendental equation (also investigated 
by Euler in 1783 \cite{e-1783}). Some of these properties will be discussed in the next section.

\section{Utility of the Lambert $W$ function}

\begin{figure}
\includegraphics[width=.9\textwidth]{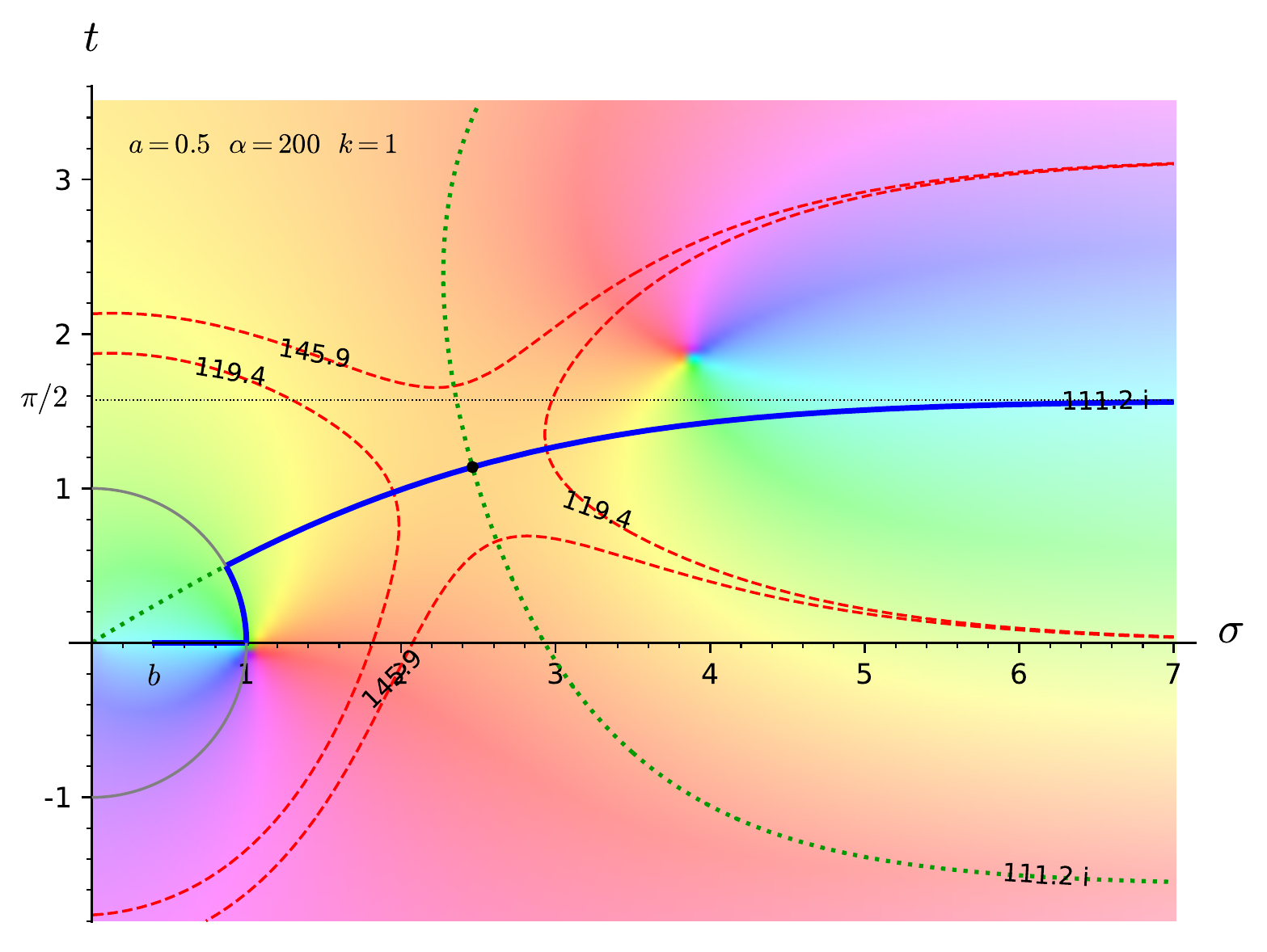}
 \caption{The function \(h_1\) for \(\alpha=200\) 
with the saddle point \(\bullet\) near \(2.46+1.14i\).
 The dashed red lines are level lines of \(\Re(h_1(\sigma+it))\)
 and the dotted green lines are the level lines where \(\Im(h_1(\sigma+it))=\Im(h_1(w_{1}(200)))\).
 The solid blue line is our path of integration. 
}\label{fig saddle 1}
\end{figure}

We first rewrite \(S_k\) such that it becomes easier to find the saddle point mentioned above. 
Recall that $f_{\alpha}(x)$ was defined as $f_{\alpha}(x) = \frac{\log^{\alpha}(x+a)}{x+a}$, 
which means that for its first two derivatives we have $f'_{\alpha}(x) =\frac{\log^{\alpha-1}(x+a)}{(x+a)^2}(\alpha-\log(x+a))$ and 
\begin{align}\label{eq ddf}
f''_{\alpha}(x) =\frac{\log^{\alpha-2}(x+a)}{(x+a)^3}\left(\alpha(\alpha-1)-3\alpha\log(a + x)+2\log^2(a + x)\right).
\end{align}
Now, since \(S_k :=\int_1^\infty  e^{2\pi ikx} f_{\alpha}''(x) dx\), with the change of variables \(y=\log(x+a)\) 
and \(b=\log(1+a)\) we can write 
\begin{align*}
S_k
&=
\int_1^\infty 
e^{2\pi ikx}
\frac{\log^{\alpha-2}(x+a)}{(x+a)^3}\left(\alpha(\alpha-1)-3\alpha\log(a + x)+2\log^2(a + x)\right) \; dx \\
&= 
\int_b^\infty e^{2\pi ik(e^y-a)}
\frac{y^{\alpha-2}}{e^{3y}}\left(\alpha(\alpha-1)-3\alpha\cdot y+2\cdot y^2\right) e^y \; dy \\
&=
\int_b^\infty e^{2\pi ik(e^y-a)+\alpha\log y}e^{-2y}
\left( \frac{\alpha(\alpha-1)-3\alpha\cdot y+2\cdot y^2}{y^2} \right)  dy.
\end{align*}
Let $h_k(y)=2\pi i k(e^y-a) + \alpha \log y$ and \(q(y)=\frac{\alpha(\alpha-1)-3\alpha\cdot y+2\cdot y^2}{y^2}\).  Then 
\begin{align}
S_k=\int_b^\infty e^{h_k(y)}e^{-2y}q(y) \; dy.
\end{align}

The function \(h_k(y)\) defined above  has a saddle point where $h_k'(y_1)=2\pi i k e^{y_1}+\alpha/y_1=0$ and the Lambert $W$ function tells us that this happens at
\(y_1={W\left(\frac{\alpha i}{2\pi k}\right)}\).
We use the principal branch $W_0$ of the Lambert $W$ function and set 
\[
w_k(\alpha)=W_0\left(\frac{\alpha i}{2\pi k}\right).
\]
We make a couple of observations concerning \(W_0(it)\) that will be useful later.
\begin{lemma}\label{lem:T}
Let \(W_0\) be the principal branch of the Lambert $W$ function.  For \(t\in(0,\infty)\), the inverse of \(I(t):=\Im(W_0(it))\) is the function $T(y)$, where for $y\in(0,\pi/2)$:
\begin{equation}\label{eq T}
T(y) = \frac{y}{\cos y} \cdot e^{y\cdot\tan(y)}.
\end{equation}
\end{lemma}

\begin{proof}
Considering the real part of $W_0(it)\cdot e^{W_0(it)}  = it$ 
we get 
\(
\Re(W_0(it))=I(t)\cdot\tan(I(t))
\), so that 
\(
W_0(it)=I(t)(\tan(I(t))+i).
\)
Using this in $W_0(it)\cdot e^{W_0(it)}  = it$, 
and considering only the imaginary parts, we obtain
\(
I(t)\cdot e^{I(t)\tan(I(t))}\frac{1}{\cos(I(t))}=t
\).
This shows that, for \(y\in [0,\pi/2)\), if we set  $T(y)=\frac{y}{\cos y} \cdot e^{y\cdot\tan(y)}$, then \(T\) is the inverse of \(I\). 
\end{proof}

It follows immediately from Lemma \ref{lem:T}  that \(T(0)=0\) and \(\lim_{y\to\pi/2} T(y)=\infty\) and that
\[
T'(y)=\left(((\tan^2(y) + 1)\cdot y + \tan(y))\cdot y+ y\tan y + 1
\right)\frac{e^{y\tan(y)}}{\cos(y)}.
\]
Hence \(T'(y)>0\) for \(y\in [0,\pi/2)\), and therefore \(I(t)>0\), for \(t\in(0,\infty]\).
Thus for \(t>0\) we have
\begin{equation}\label{eq:I lt}
0<I(t)<\frac{\pi}{2}
\end{equation}
and \(\lim_{t\to\infty}I(t)=\pi/2\).
This implies
\begin{equation}\label{eq:R gt}
\Re(W_0(it))=I(t)\cdot\tan(I(t))>0.
\end{equation}
Taking the logarithm of $W_0(it)\cdot e^{W_0(it)}  = it$ and only considering real parts we obtain
$\Re(W_0(it)) = \log(t) - \Re(\log W_0(it))$. Thus, for  \(t>1.97\) where \(|W_0(it)|>1\),
\begin{equation}\label{eq:W lt log}
\Re(W_0(it))<\log(t).
\end{equation}
We will use the following two lemmas to show that we can set \(S^*=|S_1|\).

\begin{lemma}\label{lem:deriv} 
For \(t>0\) we have \(\frac{d}{dt}\left(\Re\left(\log(W_0(it))-\frac{1}{W_0(it)}\right)\right)>0\). 
\end{lemma}

\begin{proof}
With \(W_0'(x)=\frac{W_0(x)}{x\cdot(1+W_0(x))}\) we get
\(
\frac{d}{dt}\Re\left(\log(W_0(it)-\frac{1}{W_0(it)}\right)= \frac{1}{t}\frac{\Re W_0(it)}{|W_0(it)|}>0
\),
as desired.
\end{proof}
\begin{lemma}\label{lem sin bound}
For \(k\in\N\) and \(\alpha\in[2\pi,\infty)\) we have
\(
k\cdot\sin(\Im(w_k(\alpha)))\ge 1.
\)
\end{lemma}
\begin{proof}
Representing the cosine and tangent functions in (\ref{eq T}) by the sine function, and using that  
\(\sin^{-1}x\le \frac{\pi}{3}x\), for \(x\in\left[0,\frac{1}{2}\right]\), we deduce: 
\begin{align*}
T\left(\sin^{-1}\frac{1}{2k}\right) 
\; &= \; \; \frac{\sin^{-1}\frac{1}{2k}}{\cos \left(\sin^{-1}\frac{1}{2k}\right)} 
\cdot e^{\left(\sin^{-1}\frac{1}{2k}\right)\cdot\tan\left(\sin^{-1}\frac{1}{2k}\right)}\\
&
= \; \; \frac{\sin^{-1}\frac{1}{2k}}{\sqrt{1-\sin^2 \left(\sin^{-1}\frac{1}{2k}\right)}} \cdot e^{\left(\sin^{-1}\frac{1}{2k}\right)\cdot \frac{\sin\left(\sin^{-1}\frac{1}{2k}\right)} {\sqrt{{1-\sin^2\left(\sin^{-1}\frac{1}{2k}\right)}}}} \\
& \le \; \; \frac{\frac{\pi}{6k}}{\sqrt{1-\frac{1}{4k^2}}}              \cdot e^{\frac{\pi}{6k}\cdot \frac{\frac{1}{2k}} {\sqrt{{1-\frac{1}{4k^2}}}}} 
\; \; \le \; \; \frac{\pi}{3\sqrt{4k^2-1}}              \cdot e^{\frac{\pi}{6k \sqrt{{4k^2-1}}}}\\
& < \; \; \frac{\pi}{3\sqrt{3}k}              \cdot e^{\frac{\pi}{6\sqrt{3}k^2}}
\; \; < \; \; 0.61 \cdot \frac{1}{k}  \cdot 1.4
\; \;  < \; \; \frac{1}{k}.   
\end{align*}

Applying the functions \(I\) and sine to both sides of this inequality yields 
\[
\frac{1}{2k}\le \sin\left(I\left(\frac{1}{k}\right)\right).
\]
Because of the monotonicity of \(I\) for \(\alpha\ge 2\pi\), we get
\(
\frac{1}{2}\le k\sin\left(I\left(\frac{\alpha}{2\pi k}\right)\right)=k \cdot \sin\left(\Im(w_k(\alpha))\right).
\)
\end{proof}

\begin{figure}
 \includegraphics[width=.9\textwidth]{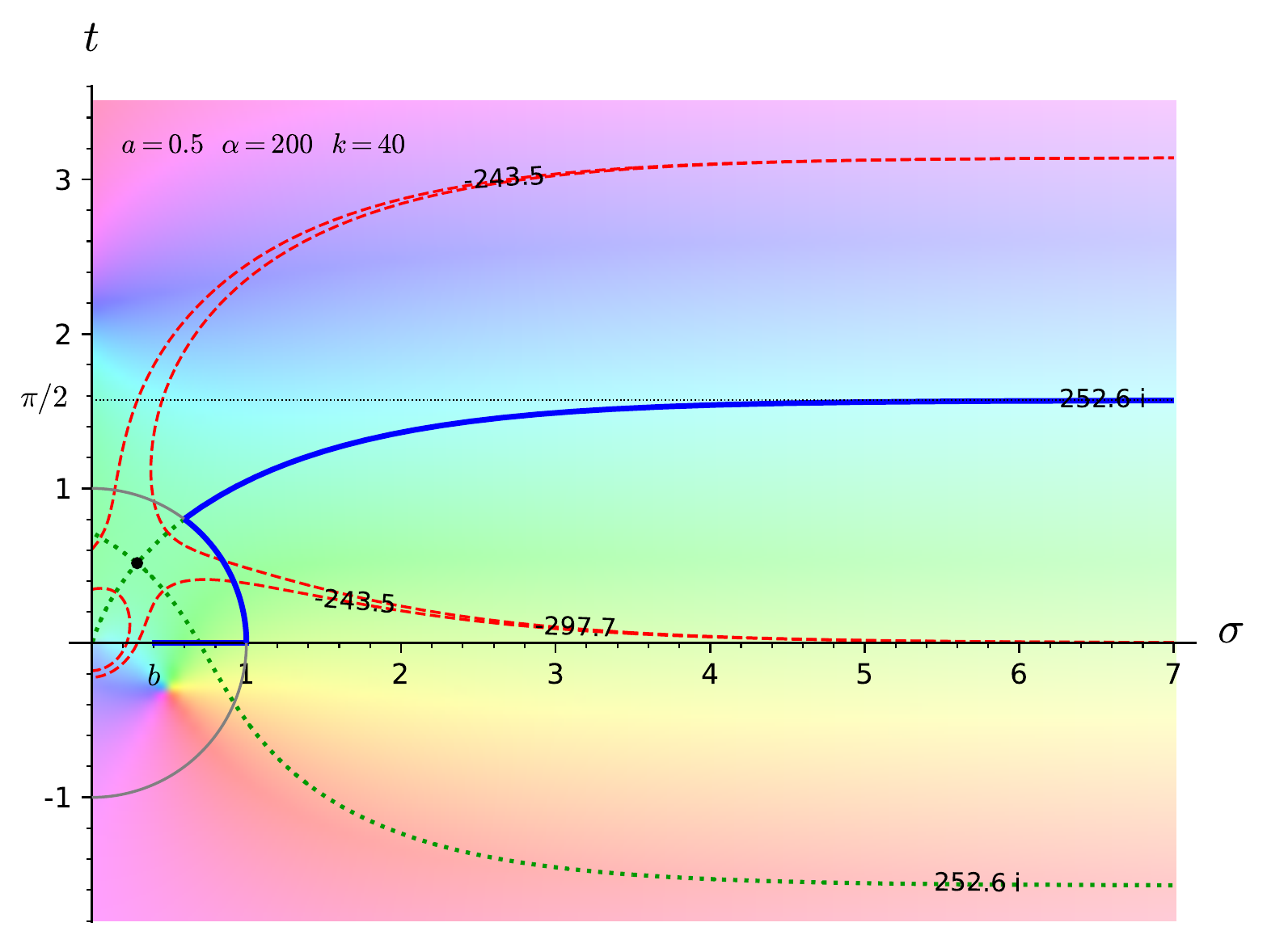}
 \caption{The function \(h_{40}\) for \(\alpha=200\)
 with the saddle point \(\bullet\) near \(0.29+0.52i\).
The dashed red lines are level lines of \(\Re(h_{40}(\sigma+it))\)
 and the dotted green lines are the level lines where \(\Im(h_{40}(\sigma+it))=\Im(h_{40}(w_{40}(200)))\).
 The solid blue line is our path of integration. 
 }\label{fig saddle 40}
\end{figure}

\section{Bounding the Integrals} 

First, recall the definitions of the quantities $|S_k|$ we are interested in: 
\[
|S_k| :=\left| \int_b^\infty e^{h_k(y)}e^{-2y}q(y) \; dy  \right|,
\]
where (as in Lemma 1) $h_k(y)=2\pi i k(e^y-a) + \alpha \log y$, \(q(y)=\frac{\alpha(\alpha-1)-3\alpha\cdot y+2\cdot y^2}{y^2}\), and \(b=\log(1+a)\).

We evaluate \(S_k\) by integrating along the contour 
that starts at \(b\), goes along the positive real axis to 1, follows the unit circle until the
level line \(\Im(h_k(y))=\Im(h_k(w_k(\alpha))\) is reached,
crosses the saddle at \(w_k(\alpha)\) (provided it is not inside the unit circle)
and continues on the level line.
From our observations (\ref{eq:I lt}) and (\ref{eq:R gt}) we know that \(0<\Im(w_k(\alpha))<\pi/2\)
and \(\Re(w_k(\alpha))>0\).
Setting \(\Im(h_k(y))=\Im(h_k(w_k(\alpha)))\) we see that as
\(\Re(y)\to\infty\) we have \(\Im(y)\to \frac{\pi}{2}\).
On the level line, from the origin to the saddle \(\Re(h_k(y))\) strictly increases and after crossing it \(\Re(h_k(y))\) strictly decreases.
See Figures \ref{fig saddle 1} and \ref{fig saddle 40}.

Following the segments of the contour we split up the integral into four parts, \(S_k=L_1+L_2+L_3+L_4\), where
\begin{itemize}
\item \(L_1=\int_b^1 e^{h_k(y)}e^{-2y}q(y)  dy \),
\item \(L_2=\int_1^u e^{h_k(y)}e^{-2y}q(y)  dy \) along the unit circle from 1 to the point $u$ where the unit circle and the level line \(\Im(y)=\Im(h_k(w_k(\alpha))\) meet,
\item \(L_3=\int_u^v e^{h_k(y)}e^{-2y}q(y)  dy \) along the level line until the point \(v\) with \(\Re(v)= 2\log\alpha\), if \(\Re(h_k(v))>\Re(h_k(u))\), otherwise \(v=u\) and \(L_3=0\),
\item \(L_4=\int_v^\infty e^{h_k(y)}e^{-2y}q(y)  dy \)  
\end{itemize}
In the following we estimate each of these four components separately.


\begin{lemma}\label{L_1}
With $L_1$ defined as above, we have
\( 
\displaystyle|L_1 | < \alpha+3+\frac{1}{\pi}.
\)
\end{lemma}

\begin{proof}
Let $y\in\R^+$. First, let us note that 
\begin{equation} \label{eq real h}
\Re(h_k(y)) =\Re(2\pi i k(e^{y}-a) + \alpha \log (y)) =-2\pi k e^{\Re(y)}\sin\Im(y) + \alpha \log |y|
\end{equation}
and because \(|e^{-2y}|\le 1\), we can write: 
\begin{align*}
|L_1| 
& = \;  \left| \int_b^1 e^{h_k(y)}e^{-2y}q(y) \; dy \right| \; \le  \;\int_b^1 |e^{h_k(y)}| \cdot |q(y)| \; dy \\
& \le  \; \int_b^1 e^{-2\pi k e^{\Re(y)}\sin\Im(y) + \alpha \log |y|}\cdot |q(y)| \; dy \; \leq \; \int_b^1 e^{\alpha \log y}\cdot |q(y)| \; dy \\
&\le \; \int_b^1 y^\alpha \left(\frac{\alpha(\alpha-1)}{y^2}+\frac{3\alpha}{y}+2\right)  dy  \; =  
\; \left[\alpha y^{\alpha-1}+3y^{\alpha}+\frac{2y^{\alpha+1}}{\alpha+1}\right]_b^1 \\
&= \; \left(\alpha+3+\frac{2}{\alpha+1}\right)-\left( \log(a+1)^{\alpha-1}+3\log(a+1)^{\alpha}+\frac{2\log(a+1)^{\alpha+1}}{\alpha+1}\right) \; < \; \alpha+3+\frac{1}{\pi}\qedhere
\end{align*}
\end{proof}

\begin{lemma}\label{L_2}
With $L_2$ defined as above, we have
$\displaystyle |L_2 | < (\alpha^2+2\alpha+2)\frac{\pi}{2}. $
\end{lemma}

\begin{proof}
First observe that on the unit circle we have \(|y|=1\) and \(0\le \Im(y)\le 1\), and thus, by (\ref{eq real h}), 
\begin{equation}
{\Re(h_k(y))}
={-2\pi k e^{\Re(y)}\sin\Im(y) + \alpha \log |y|} 
={-2\pi k e^{\Re(y)}\sin\Im(y)} \le 0.
\end{equation}
Moreover, we can write 
\begin{equation}\label{eq:bound q}
|q(y)| \; \le \; \left|\frac{\alpha(\alpha-1)}{y^2}\right|+\left|\frac{3\alpha}{y}\right|+2
\; = \; \frac{\alpha(\alpha-1)}{|y^2|}+\frac{3\alpha}{|y|}+2 \; = \; \alpha^2+2\alpha + 2.
\end{equation}
Since replacing $u$ by $i$ can only extend the path of integration, the last bound directly gives: 
\begin{equation*}\label{eq bound C}
|L_2| \; \le \; \int_1^i \left|e^{h_k(y)}e^{-2y}q(y)\right| \; dy \; \le \; \int_1^i |e^{h_k(y)}|\cdot|e^{-2y}|\cdot|q(y)| \; dy\; = \; (\alpha^2+2\alpha+2)\frac{\pi}{2}.\qedhere
\end{equation*}
\end{proof}

\begin{lemma}\label{L_4}
With $L_4$ defined as above, we have
$|L_4 | < \alpha^2+2\alpha+2$.
\end{lemma}

\begin{proof}
The remainder of our integration path follows the level line with \(\Im(y)=\Im(w_k(\alpha))\).
Here we have $|y|\ge 1$ and thus, as in (\ref{eq:bound q}) we have $q(y)\le \alpha^2+2\alpha+2$.  We get
\begin{equation}\label{eq right}
\left| \int e^{h_k(y)}e^{-2y}q(y)dy \right|
\le \int |e^{h_k(y)-2y}|(\alpha^2+2\alpha+2) dy
= (\alpha^2+2\alpha+2) \int e^{\Re(h_k(y)-2y)} dy
\end{equation}
By (\ref{eq:W lt log}) we have \(\Re(w_k(\alpha))<\log(\alpha)\).  So \(v\) lies to the right of the saddle.
For \(\Re(y) \ge 2\log\alpha\) we have
\begin{equation}\label{eq pi and e}
\Re(h_k(y)) = -2\pi k e^{\Re(y)}\sin\Im(y) + \alpha \log |y| \le  -\pi e^{\Re(y)} + \alpha \log |y| < 0.
\end{equation}
To see this, just note that in our region $\Re(y) > 4 \log \log \Re(y)$ and thus also
$\Re(y)/2 + \log \pi > 2 \log \log \Re(y) + \log 2 > \log \log (2 \Re(y)^2) = \log \log (\Re(y)^2 + \Im(y)^2) = \log \log |y|$, due to the concavity of both the logarithmic 
and the Lambert $W$ functions.
Plugging (\ref{eq pi and e}) into (\ref{eq right}) yields:
\begin{align*}
|L_4|
&\le (\alpha^2+2\alpha+2)\int_c \left| e^{-\pi e^{\Re(y)} + \alpha \log |y|}\right|e^{\Re(-2y)} dy \le (\alpha^2+\alpha+2)\int_c e^{-2y}  dy < \alpha^2+2\alpha+2.\qedhere
\end{align*}
\end{proof}

\begin{lemma}\label{L_3}
With $L_3$ defined as above, we have
\[
|L_3 | < e^{\Re\left(-\frac{\alpha}{w_k(\alpha)}+\alpha \log w_k(\alpha)\right)} (\alpha^2+2\alpha+2)\sqrt{4\log^2\alpha+\frac{\pi^2}{4}}
\]
\end{lemma}

\begin{proof}
Here the curve \(c\) has its (real) maximum at the saddle point \(w_k(\alpha)\), where \(h_k(w_k(\alpha))=2\pi i k e^{w_k(\alpha)}+\alpha/w_k(\alpha)=0\). This allows us to bound the real part of $h_k(y)$ as:
\begin{equation}
\Re(h_k(y))
\le \Re(h_k(w_k(\alpha))
= \Re\left(2\pi i k(e^{w_(\alpha)}-a)+\alpha \log w_k(\alpha)\right)
= \Re\left(-\frac{\alpha}{w_k(\alpha)}+\alpha \log w_k(\alpha)\right),
\end{equation}
which means that we can estimate the integral as: 
\begin{align*}
|L_3|
&=\; \left| \int_u^v e^{h_k(y)}e^{-2y}q(y) \; dy \right| \le e^{\Re\left(-\frac{\alpha}{w_k(\alpha)}+\alpha \log w_k(\alpha)\right)} (\alpha^2+\alpha+2) \int_u^v dy\\
&\le \; e^{\Re\left(-\frac{\alpha}{w_k(\alpha)}+\alpha \log w_k(\alpha)\right)} (\alpha^2+2\alpha+2) \int_0^{2\log\alpha+\pi/2} dy \\
&= \; e^{\Re\left(-\frac{\alpha}{w_k(\alpha)}+\alpha \log w_k(\alpha)\right)} (\alpha^2+2\alpha+2)\sqrt{4\log^2\alpha+\frac{\pi^2}{4}}.\qedhere
\end{align*}
\end{proof}

Putting these four Lemmas together, we immediately get the following bound:
\begin{align}
 |S_k| &\le \; |L_1| + |L_2| + |L_3| + |L_4| \nonumber\\
&< \; \alpha+3+\frac{1}{\pi}+(\alpha^2+2\alpha+2)\left(1+\frac{\pi}{2}+ e^{\Re\left(-\frac{\alpha}{w_k(\alpha)}+\alpha \log w_k(\alpha)\right)}
\sqrt{4\log^2\alpha+\frac{\pi^2}{4}}. 
\right)
\label{eqn:S bound}
\end{align}

\section{The Final Bound}

Now we can prove the following general result: 

\begin{theorem}\label{theo:bound}
For \(\alpha\ge 2\pi\) and \(a\in(0,1]\)
let us denote by \(\gamma_\alpha(a)\) the fractional Stieltjes constants and write 
\(C_\alpha(a)=\gamma_\alpha(a)-\frac{\log^\alpha a}{a}\).  
If we set $w(\alpha):={W_0\left(\frac{\alpha i}{2\pi}\right)}$, where
\(W_0\) is the principal branch of the Lambert $W$ function, then
\[
|C_\alpha(a)| \; < \; \alpha^2+ \frac{3}{4} \alpha^2 \log \alpha \cdot \left|e^{\alpha(\log w(\alpha)-1/w(\alpha))}\right|. 
\]

\end{theorem}
\begin{proof}
From Lemmas \ref{lem:deriv} and \ref{L_4} it follows that the quantities \(S_k\) decrease as \(k\) increases.
Thus we can set \(S^*=|S_1|\) in (\ref{eq:S}), and with the help of the bound (\ref{eqn:S bound})  we can rewrite (\ref{eq:Calpha}) as:
\begin{align*}
|C_\alpha(a)|   
& <   \log^\alpha(1+a)\left|\frac{1}{2}\frac{1}{1+a} - \frac{\log(1+a)}{\alpha+1} + \frac{1}{12(1+a)^2}\right|
+\left| \frac{1}{12}\frac{\log^{\alpha-1}(1+a)}{(1+a)^2}\alpha\right|  +  \frac{1}{6}|S_1| \\
& <   1+\frac{1}{12}\alpha  + \frac{1}{6}\left[\alpha+3.3+(\alpha^2+2\alpha+2)\left( \left(1+\frac{\pi}{2} \right)+ \left(2 \log \alpha+ \frac{\pi}{2} \right) \cdot e^{\Re\left(-\frac{\alpha}{w_1(\alpha)}+\alpha \log w_1(\alpha)\right)}\right)\right] \\
& < 1.55 + \frac{\alpha}{4} + \frac{1}{2}\left( \alpha^2+2\alpha+2 \right)  + \frac{1}{6} \left(\frac{3}{2}\alpha^2 \right)(3 \log \alpha) \cdot e^{\Re\left(-\frac{\alpha}{w_1(\alpha)}+\alpha \log w_1(\alpha)\right)}\\
&  < \; \alpha^2+ \frac{3}{4} \alpha^2 \log \alpha \cdot \left|e^{\alpha(\log w(\alpha)-1/w(\alpha))}\right|,
\end{align*}
since $\frac{1}{6}(1 + \frac{\pi}{2}) < \frac{1}{2}$ and in our range we have $2.55 + \frac{5}{4}\alpha  < \frac{1}{2}\alpha^2$. This establishes the bound. 
\end{proof}

Note that the main term of the bound in Theorem \ref{theo:bound} differs only by a factor of 
\(\alpha^2\log\alpha\) from the conjectured bound given in \cite{fpsbound}:
\begin{equation}\label{eq:conjecture}
|C_\alpha(a)|\le 2\left|\displaystyle e^{\alpha\left(\log w(\alpha)-{1}/{w(\alpha)}\right)}\right|.
\end{equation}

\begin{figure}
 \includegraphics[width=.9\textwidth]{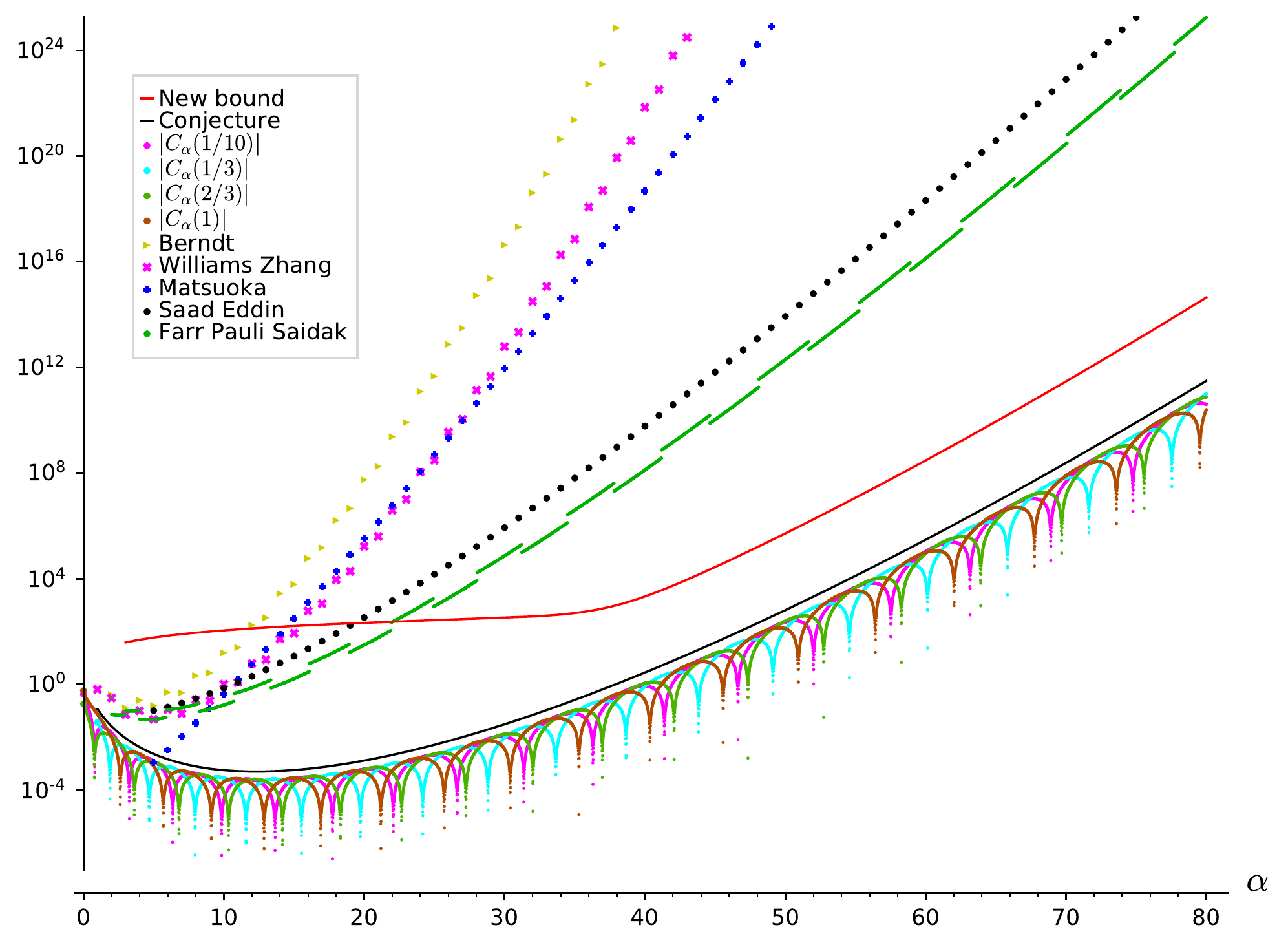}
 \caption{On a logarithmic scale we show the absolute values of the Stieltjes constants \(\gamma_\alpha=C_\alpha(1)\),
\(C_\alpha(2/3)\),
\(C_\alpha(1/3)\) and 
\(C_\alpha(1/10)\) 
along with the bounds by Berndt \cite{berndt1972}, Williams and Zhang \cite{wz1}, Matsuoka  \cite{matsuoka} and Saad Eddin \cite{saad}, 
our bound and conjecture from \cite{fpsbound} as well as our new bound from Theorem \ref{theo:bound}.
}\label{fig bounds}
\end{figure}

In Figure \ref{fig bounds} we compare Theorem \ref{theo:bound} and (\ref{eq:conjecture}) with previously known bounds for \(\gamma_\alpha=C_\alpha(1)\).
For \(m\in\N\) we have:
\begin{enumerate}
\item the bound by Berndt \cite{berndt1972}:
\(
|\gamma_m|\le  \frac{(3+(-1)^{m})(m-1)!}{\pi^m}
\)
\item the bound by Williams and Zhang \cite{wz1}:
\(
|\gamma_m|\le  \frac{(3+(-1)^{m})(2m)!}{m^{m+1}(2\pi)^m}
\)
\item the bound by Matsuoka \cite{matsuoka} which holds for $m>4$:
\(
|\gamma_m|<10^{-4}(\log m)^m
\)
\item the bound by Saad Eddin \cite{saad}:
Let
$\theta(m)=\frac{m+1}{\log\frac{2(m+1)}{\pi}}-1$ then
\[
|\gamma_m|
\le m!\cdot2\sqrt{2}e^{-(n+1)\log\theta(m)+\theta(m)
\left(\log\theta(m)+\log\frac{2}{\pi e}\right)}
\left(1+2^{-\theta(m)-1}\frac{\theta(m)+1}{\theta(m)-1}\right).
\]
\item our bound from \cite{fpsbound}:  
For \(\alpha\in(0,\infty)\) let \(x = \frac{\pi}{2} e^{W_0\left( \frac{2(\alpha+1)}{\pi} \right)}\) 
then
\[
|\gamma_\alpha| \leq  \frac{(3+(-1)^{n+1})
\Gamma(\alpha+1)}{(2\pi)^{n+1}(n+1)^{\alpha+1}} \frac{(2(n+1))!}{(n+1)!}
\mbox{ where }
n=\left\{
\begin{array}{ll}
\lfloor x \rceil & \mbox{if } x<\alpha\\
\lceil\alpha-1\rceil & \mbox{otherwise}
\end{array}
\right.
\] 
\end{enumerate}
\begin{remark*}
With (\ref{eq:W lt log}) we get \(\Re(\log(w(\alpha))-1/w(\alpha))<\Re(\log w(\alpha))<\log\log(\alpha)\).
Hence
\[
\left|\displaystyle e^{\alpha\left(\log w(\alpha)-{1}/{w(\alpha)}\right)}\right|
=e^{\alpha\Re\left(\log w(\alpha)-{1}/{w(\alpha)}\right)}
<(\log\alpha)^\alpha
\]
Thus the main term of Matsuoka's bound 
follows from (\ref{eq:conjecture}).
\end{remark*}
\section{Acknowledgments}
We'd like to thank the anonymous referee for several useful comments and for pointing out a couple of miscalculations.
All plots were created with the computer algebra system SageMath \cite{sage}.


\end{document}